\documentclass[12pt]{amsart}

\usepackage {amsmath}
\usepackage[applemac]{inputenc}
\usepackage {amsfonts, amssymb}
\usepackage {latexsym}
\usepackage {harvard}

\newtheorem {theorem} {Theorem} [section]
\newtheorem {lemma} {Lemma} [section]
\newtheorem {prop} {Proposition} [section]

\newtheorem{preremark}{Remark}[section]
\newenvironment{remark}%
  {\begin{preremark}\rm}{\end{preremark}}
   \newtheorem{preremark1}{Example}[section]
\newenvironment{example}
  {\begin{preremark1}\rm}{\end{preremark1}}
  \newtheorem{preremark2}{Definition}[section]
\newenvironment{defn}
  {\begin{preremark2}\rm}{\end{preremark2}}

\newcommand {\R} {\mathbb{R}}
\def\be{\begin{equation}}
\def\ee{\end{equation}}

\begin{document}

\title[Randomly flashing ratchet]{Analytical aspects of the Brownian motor effect in randomly flashing ratchets}

\author[D.Vorotnikov]{Dmitry Vorotnikov}

\address{CMUC, Department of
Mathematics, University of Coimbra, 3001-454 Coimbra, Portugal}{}
\email{mitvorot@mat.uc.pt}

\thanks{The research was partially supported by CMUC and FCT (Portugal), through
European program COMPETE/FEDER under the project PEst-C/MAT/UI0324/2011. The author thanks David Kinderlehrer who introduced him to ratchets, and Philippe Verkerk for a discussion on the topic. He also thanks the unknown referees for detailed comments which helped to significantly improve this article.}




\begin{abstract} The muscle contraction, operation of ATP synthase, maintaining the shape of a cell are believed to be secured by motor proteins, which can be modelled using the Brownian ratchet mechanism. We consider the randomly flashing ratchet model of a Brownian motor, where the particles can be in two states, only one of which is sensitive the applied spatially periodic potential (the mathematical setting is a pair of weakly coupled reaction-diffusion and Fokker-Planck equations). We prove that this mechanism indeed generates unidirectional transport by showing that the amount of mass in the wells of the potential decreases/increases from left to right. The direction of transport is unambiguously determined by the location of each minimum of the potential with respect to the so-called diffusive mean of its adjacent maxima. The transport can be generated not only by an asymmetric potential, but also by a symmetric potential and asymmetric transition rates, and as a consequence of the general result we derive explicit conditions when the latter happens. When the transitions are localized on narrow active sites in the protein conformation space, we find a more explicit characterization of the bulk transport direction, and infer that some common preconditions of the motor effect are redundant.

\end{abstract}

\maketitle

\section {Introduction}

\emph{Brownian motors} are nano-scale or molecular machines which can produce directed motion when the average force and average temperature gradient are zero \cite{ast,bm}. At first glance, the existence of such devices seems to be a paradox, and discretization of the idea really did lead to Parrondo's paradox in game theory \cite{ameng,har}. Typically, the mechanism involved (\emph{ratchet}) is based on an interplay between the Brownian motion (\emph{diffusion}), an asymmetric (ratchet-like) \emph{potential}, and \emph{nonequilibrium} of the system due to chemical or thermal \emph{fluctuations}.

The ratchet principle is ubiquitous and appears everywhere from political system to famine cycles, from production strategy to cultural studies. Motor proteins, which provide muscle contraction (myosin, kinesin, dynein), ATP synthase, as well as membrane-bound motor proteins maintaining the shape of a cell, can be modelled using the Brownian ratchet mechanism \cite{modbio,jul,mod2,mod3,mod1}. A ribosome can also be considered as a Brownian ratchet device \cite{sp}.

Motor proteins can attach to and detach from a substrate
of vectorial symmetry \cite{modbio,mod1} under the action
of a chemical energy source. This leads to the following boundary value problem \cite{modbio,jul,modmain1,modmain}:
\be \label{eq1}
\left\{\begin{array}{ll} p_t-\sigma p_{xx}- \kappa(\psi_x p)_x=  \nu P - \eta p, & x\in (0,1),\ t\geq 0, \\
P_t-\varsigma P_{xx}= -\nu P + \eta p,
& x\in (0,1), \ t\geq 0,\\
\sigma p_{x}+ \kappa\psi_x p=0, & x=0,1, \ t\geq 0,\\
P_{x}=0, & x=0,1, \ t\geq 0,\\
p\geq 0, \
P\geq 0, \\ \int\limits_0^1 p(t,x) + P(t,x)\, dx =1, \ t\geq 0.
\end{array}\right.
\ee
Here $p(t,x)$ and $P(t,x)$ are the unknown densities of the particles in ``attached'' and ``detached'' states, resp., at a time $t$ and a spatial point $x$; $\sigma>0$ and $\varsigma>0$ are the diffusion coefficients of ``attached'' and ``detached'' particles, resp.; $\psi(x)$ is the potential; $\kappa>0$ is a coefficient inversely proportional to temperature; $\nu(x)>0$ and $\eta(x)>0$ are the rates of transition from one state to another, i.e. $\nu$ indicates the probability of seizing a motor protein by a ``detached'' particle (located at a spatial point $x$), and $\eta$ expresses the probability of losing its motor protein for an ``attached'' particle. A typical ``ratchet-like'' potential $\psi$ with $k$ teeth, $k>1$, is $1/k$-periodic in $x$ and has a unique local (and, hence, global) minimum within each period.

Let us compare this model with a somehow simpler one called the \emph{flashing ratchet} \cite{ast}. Here there is only one state, $\rho$ is the unknown density of particles, $\sigma$ is the diffusion coefficient, and $\psi$ is the potential,
which is switched on and off cyclically:
\be \label{fr1}
\left\{\begin{array}{ll} \rho_t= \sigma\rho_{xx}+ h(t)(\psi_x\rho)_x,
\qquad x\in (0,1),\\
\sigma\rho_{x}+ h(t)\psi_x\rho=0, \qquad x= 0,1,\\
\rho\geq 0,\\ \int\limits_0^1 \rho(x,t) dx =1, \\ h(t)= 1,\ nT<t\leq nT+T_{tr},\ n=0,1,\dots,\\ h(t)=0, \ nT+T_{tr}<t\leq nT+T,\ n=0,1,\dots.
\end{array}\right.
\ee
Here each particle is potential-sensitive and potential-insensitive for a priori known moments of time, whereas in the first model this is determined by random attachments and detachments of the motor protein. Therefore model \eqref{eq1} is sometimes referred to as the \emph{randomly flashing ratchet}.

Another model related to \eqref{eq1} is the \emph{collaborative ratchet} \cite{chipm,jul,mod2,mod3,mod1} where the particles in two states are sensitive to two different potentials $\psi$ and $\Psi$ which help each other to achieve the motor effect:
\be \label{eq1two}
\left\{\begin{array}{ll} p_t-\sigma p_{xx}- \kappa(\psi_x p)_x=  \nu P - \eta p, & x\in (0,1),\ t\geq 0, \\
P_t-\varsigma P_{xx}- \kappa(\Psi_x P)_x= -\nu P + \eta p,
& x\in (0,1), \ t\geq 0,\\
\sigma p_{x}+ \kappa\psi_x p=0, & x=0,1, \ t\geq 0,\\
\varsigma P_{x}+ \kappa\Psi_x P=0, & x=0,1, \ t\geq 0,\\
p\geq 0, \
P\geq 0, \\ \int\limits_0^1 p(t,x) + P(t,x)\, dx =1, \ t\geq 0.
\end{array}\right.
\ee
This model is relevant not only in connection with biology, but also in transport of cold rubidium atoms \cite{cold}.

A diffusion-free problem somehow similar to \eqref{eq1two} is a model of transport of neurofilaments in axons; its behaviour near equilibrium was investigated by \citeasnoun{af1}.

The mathematical studies of Brownian ratchet models start with the question whether they indeed generate unidirectional transport (which is observed in experiments and simulations), and what assumptions are needed for that. For example, model \eqref{eq1} is too general to produce transport with any choice of parameters: setting $$\eta/\nu=\exp{(\kappa \psi/\sigma)},$$ we observe that the total amount of mass in the ``wells'' of the potential, i.e. in the segments $[\frac i k,\frac {i+1} k]$, $i=0,\dots,k-1$, eventually with the course of time tends to $1/k$.
This shows that the ratchets should be ``tuned'' to work well.

Comparing eventual distribution of mass between the ``wells'' of the potential $\psi$, it is possible to mathematically endorse the motor effect \cite{kind1}. A left-to-right chain of inequalities in this distribution would mean unidirectional transport.

The occurrence of transport in model \eqref{eq1two} for small $\sigma=\varsigma$ and certain interplay between the potentials was shown by \citeasnoun{chip}. Some of the results of that paper --- in particular, existence of unique solutions to \eqref{eq1two} and to the stationary problem
\be \label{eqbtwo1}
\left\{\begin{array}{ll} -\sigma p_{xx}- \kappa(\psi_x p)_x= \nu P - \eta p, & x\in (0,1), \\
-\varsigma P_{xx}- \kappa(\Psi_x P)_x= -\nu P + \eta p,
& x\in (0,1), \\
\sigma p_{x}+ \kappa\psi_x p=0, & x=0,1, \\
\varsigma P_{x}+ \kappa\Psi_x P=0, & x=0,1, \\
p\geq 0, \
P\geq 0, \\ \int\limits_0^1 p(x) + P(x)\, dx =1,
\end{array}\right.
\ee
and eventual convergence of solutions to \eqref{eq1two} to the solutions of \eqref{eqbtwo1} --- are valid for $\Psi\equiv\, const$ as well. The transport result was generalized to the multi-state systems with several interacting non-flat potentials by \citeasnoun{hast} and by \citeasnoun{pet}.

The analytical proof of the motor effect for the flashing ratchet \eqref{fr1} was given by \citeasnoun{vorr}, based on a framework developed by \citeasnoun{kind1}. The unidirectional transport occurs when the potential is asymmetric, and its direction is determined by the location of the minima of the potential  with respect to the centres of the corresponding wells. A homogenization approach to the flashing ratchet \eqref{fr1} was proposed by \citeasnoun{petfg}. This approach was applied to the randomly flashing ratchet \eqref{eq1} by \citeasnoun{peth} \citeaffixed{mer}{see also}. There were presented examples of transport in the case of a sawtooth potential (which is asymmetric) and constant transition rates (which are obviously symmetric), and with a symmetric potential and asymmetric transition rates. We are not aware of any works with rigorous mathematical evidences of transport for the non-homogenized problem \eqref{eq1}.

In this paper, we show that the motor effect in model \eqref{eq1} is due to a subtle interplay between the asymmetries of the potential $\psi$ and the transition rate $\nu$. More precisely, an asymmetry of the transition function $\nu$ yields a deviation of the so-called \emph{diffusive means} of the edge points of the wells of the potential --- these edge points are the maxima of $\psi$ --- from the centres of the wells. The direction of transport is determined by the location of the minima of the potential with respect to that biased centres (the diffusive means of the adjacent maxima). The ratchet is tuned when the influence of the potential on the particles which carry the motor protein dominates their diffusion (which is a natural assumption since these complex particles are larger and thus diffuse slower).

Among the consequences of the general result, the following observation is of special importance in connection with motor proteins. The function $\nu(x)$ is conformation-dependent, and is frequently assumed \cite{jul} to dramatically increase on narrow regions corresponding to the areas in the conformation space where the reattachment activity of motor proteins is overwhelming, and to be negligibly small everywhere else. In this particular framework, we infer that the bulk transport direction is dictated by the location of the minima of the potential with respect to that active sites, and observe that several preconditions of the motor effect found in biophysical literature are unnecessary.

Let us briefly announce the main biological insights of the paper. We highlight the main factors that manage the operation of biomotors: chemical potential, architecture of the involved proteins, structural changes in the system, including  attachment/detachment of motor proteins, and diffusion. We reveal and investigate the intricate interplay between these factors, and observe that their subtle alterations may change collective behaviour of the system, in particular, the transport direction. Previous analytical studies revolved around the primacy of asymmetry of the chemical potential, whereas the influence of the transition rates was regarded as a secondary factor  \cite{chip,hast,pet} or ignored \cite{kind1,petfg,vorr}.  Our key finding is that the spatial distribution of the active sites (where the attachment of motor proteins occurs) is as important as the form of the potential (\citeasnoun{peth} give an example in support of this claim, but here we provide a comprehensive proof). We can theoretically anticipate the behaviour (for example, the bulk transport direction) of motor proteins and their cargoes. We contribute to the understanding of the mechanism of motor proteins, and demonstrate that usefulness of qualitative mathematical analysis cannot be underrated prior to designing biomotor-powered hybrid devices and, eventually, artificial biomotors \cite{pax}.

The paper is organized as follows. In Section 2, we set the problem more rigorously, define the notion of the diffusive mean, and  formulate the main result (Theorem \ref{mainthmr}).  Section 3 is a discussion of consequences of Theorem \ref{mainthmr}. In particular, we investigate the case mentioned in the previous paragraph. In this section, we also study some properties of the diffusive mean and evaluate it. In Section 4, we introduce a semidiscretized device, the \emph{squeezing ratchet}, and prove that it generates unidirectional transport. The proof of the main theorem is provided in Section 5, where we show that the squeezing ratchet and the original randomly flashing ratchet have similar behaviour. Section 6 is a biologically-oriented discussion of the results. We opted for moving all the proofs except the one of Theorem \ref{mainthmr} to several appendices.

\section {Preliminaries}

We consider the stationary boundary value problem for the randomly flashing ratchet equation with Neumann boundary conditions
\be \label{eqb1}
\left\{\begin{array}{ll} -\sigma p_{xx}- \kappa(\psi_x p)_x= \nu P - \eta p, & x\in (0,1), \\
-\varsigma P_{xx}= -\nu P + \eta p,
& x\in (0,1), \\
\sigma p_{x}+ \kappa\psi_x p=0, & x=0,1, \\
P_{x}=0, & x=0,1, \\
p\geq 0, \
P\geq 0, \\ \int\limits_0^1 p(x) + P(x)\, dx =1,
\end{array}\right.
\ee which describes the eventual distribution of particles subjected to the action of the ratchet.

The potential $\psi(x)$ and the transition rates $\nu(x)>0$ and $\eta(x)>0$ are assumed to be smooth scalar functions on $[0,1]$ of period $1/k$, with $k>1$ being a fixed integer. The potential $\psi$ should have maxima at points $x_i$ and minima at points $a_i$, and be monotonic (without zero slopes) between these points, where \be x_i=\frac {i-1} k,\ i=1,\dots, k+1,\ee \be a_i= a+x_i,\ i=1,\dots, k.\ee The positive parameter $a$ should be less than $1/k$.

We use the bra-ket notation $\langle \mu,f\rangle=\int\limits_A^B f \, d\mu$, where
 $\mu\in C^*[A,B]$, and $f$ is a continuous function on $[A,B]$.

 The symbol $\delta_x$ denotes the Dirac delta centered at $x\in\R$. The symbol $C$ will stand for a generic positive constant that can take different values in different lines.

\begin{defn} Unless otherwise specified, we say that $u\in L_1(A,B)$ is a solution to the problem \be\notag\left\{\begin{array}{ll}
(k_1 u)_{xx}(x)+
(k_2 u)_x(x)+k_3(x)u(x)= \Theta(u)(x),
& x\in (A,B), \\
(k_1 u)_{x}(A)=(k_1 u)_{x}(B)=0,
\end{array}\right. \ee
where functions $k_1,k_2,k_3\in C[A,B]$, $k_2(A)=k_2(B)=0$, and a linear operator $\Theta: L_1(A,B)\to C^*[A,B]$
are prescribed, if \begin{multline} \label{ewu}
\int\limits_A^B k_1(x)u(x)\varphi_{xx}(x)-k_2(x)u(x)\varphi_x(x)+k_3(x)u(x)\varphi(x) \,dx \\ =
\langle \Theta(u), \varphi\rangle
\end{multline} for any $\varphi\in C^2[A,B]$, $\varphi_x(A)=\varphi_x(B)=0$.
\end{defn}

Note that \eqref{ewu} already includes the Neumann boundary condition.

Existence of a unique solution $(p,P)\in W^2_1(0,1)\times W^2_1(0,1)$ to \eqref{eqb1} follows from \possessivecite{chip} results  (their paper considers the case $\sigma=\varsigma$, but a change of variables yields the general case).

\begin{defn} Let $A$ and $B$ be real numbers, and $\phi(x)> 0$ be a continuous scalar function on $[A,B]$. The number $s\in (A,B)$ is called the \emph{$\phi$-diffusive mean} of $A$ and $B$ provided the system \be \label{defm}
\left\{\begin{array}{ll}
\phi U - U_{xx}=  \delta_s,
& x\in (A,B), \\
U(A)=U(B), \\
U_{x}(A)=U_{x}(B)=0
\end{array}\right.
\ee has a solution $U$.  \end{defn}

\begin{prop} \label{pr1} The $\phi$-diffusive mean always exists and is unique. \end{prop}

For any integrable scalar function $\Phi$ on $(0,1)$, we denote \be\widehat \Phi _i=\int\limits_{x_i}^{x_{i+1}}\Phi (x)\,dx,\ i=1,\dots, k.\ee

The main result of the paper is

\begin{theorem} \label{mainthmr} Let $S$ be the $\nu/\varsigma$-diffusive mean of $x_1=0$ and $x_2=1/k$.
If $a<S$, then, for sufficiently small $\sigma$ and sufficiently large $\kappa$,  \be\label{mainr}\widehat p _1>\widehat p _2>\dots>\widehat p_k,\ee \be\label{mainr1} \widehat P _1>\widehat P _2>\dots>\widehat P_k.\ee
\end{theorem}

\begin{remark} We will even prove that \be \label{qinm} P(x)> P(x+1/k), \ 0\leq x\leq 1-1/k, \ee which is stronger than \eqref{mainr1}.  \end{remark}

We can renormalize \eqref{eqb1} to get
\be \label{eqb2}
\left\{\begin{array}{ll} -\sigma p_{xx}- \kappa(\psi_x p)_x= \nu P - \eta p, & x\in (0,1), \\
-\varsigma P_{xx}= -\nu P + \eta p,
& x\in (0,1), \\
\sigma p_{x}+ \kappa\psi_x p=0, & x=0,1, \\
P_{x}=0, & x=0,1, \\
p\geq 0, \
P\geq 0, \\ \int\limits_0^1 \eta(x)p(x) + \nu(x)P(x)\, dx =2.
\end{array}\right.
\ee
Integration of \eqref{eqb2} implies
\be\int\limits_0^1 \eta(x)p(x) - \nu(x)P(x)\, dx =0.\ee
Since $0$ and $1$ are maxima of $\psi$, we have $\psi_x(0)=\psi_x(1)=0$.
Thus, \eqref{eqb2} is equivalent to
\be \label{eqb3}
\left\{\begin{array}{ll} -\sigma p_{xx}- \kappa(\psi_x p)_x= \nu P - \eta p, & x\in (0,1), \\
-\varsigma P_{xx}= -\nu P + \eta p,
& x\in (0,1), \\
p_{x}=0, & x=0,1, \\
P_{x}=0, & x=0,1, \\
p\geq 0, \
P\geq 0, \\ \int\limits_0^1 \eta(x)p(x) \, dx=1, \\ \int\limits_0^1 \nu(x)P(x)\, dx =1.
\end{array}\right.
\ee

Since problem \eqref{eqb1} is linear, it is enough to prove Theorem \ref{mainthmr} for the renormalized problem \eqref{eqb3}.

\section{Properties of diffusive means and applications of Theorem 2.1}

Theorem \ref{mainthmr} implies that, if the diffusion of potential-sensitive particles is slow, and the influence of the potential is strong (or the temperature is low), and if $a<S$, then, given any initial distribution of density, the mass of particles of each kind in the wells will eventually decrease from left to right, i.e. the motor effect is present. In order to present more specific consequences of Theorem \ref{mainthmr}, we need to better understand the notion of the diffusive mean.

\begin{example} \label{ex1} Let $$\tilde\phi(x)=\phi(A+B-x).$$ Let $s$ be the $\phi$-diffusive mean of $A$ and $B$, and let $\tilde s$ be the $\tilde\phi$-diffusive mean of $A$ and $B$. Then \be\label{til} s= A+B-\tilde s.\ee
In particular, if $\phi$ is a constant function or merely \be\phi(x)= \tilde\phi(x), \ A\leq x\leq \frac {A+B} 2,\ee then, due to \eqref{til}, the diffusive mean $s$ coincides with the arithmetic mean $\frac {A+B} 2$. In other words, $s$ is located in the centre of the segment $[A,B]$.\end{example}

In general position, the diffusive mean is distorted from the centre of the segment. We now evaluate the direction of the bias.

\begin{prop} \label{ex2} Let \be\phi(x)\leq \tilde\phi(x), \ A\leq x\leq \frac {A+B} 2,\ee and the inequality is strict at least at one point. Then the $\phi$-diffusive mean of $A$ and $B$ is strictly larger than the arithmetic mean $\frac {A+B} 2$. \end{prop}

A symmetry argument shows that if \be\phi(x)\geq \tilde\phi(x), \ A\leq x\leq \frac {A+B} 2,\ee and the inequality is strict at least at one point, then the diffusive mean is strictly less than $\frac {A+B} 2$.

Theorem \ref{mainthmr} indicates that the direction of transport generated by the ratchet is merely determined by the location of $a$ with respect to $S$ (in other words, by the direction of the vector $
\overrightarrow{Sa}$). Indeed, we have just pointed that out for $a<S$, and if $a>S$, then
$\widehat p _1<\widehat p _2<\dots<\widehat p_k$, $\widehat P _1<\widehat P _2<\dots<\widehat P_k$ --- to see this it suffices to make the change of variables $x\to 1-x$, to apply the reasoning of Example \ref{ex1} with $\phi=\frac \nu \varsigma$, $A=0$ and $B=1/k$, and to take into account that $\tilde\phi(x)=\phi(\frac 1 k  - x)=\phi(1-x)$ due to periodicity.

\begin{example} Let  the potential $\psi$ be symmetric, i.e. $\psi(x)=\psi(1-x)$, which can only happen when $a=\frac 1{2k}$. The one-state flashing ratchet cannot generate transport in this case. However, the randomly flashing ratchet can. Assume that \be\label{inex1}\nu(x)\leq \nu(-x+1/k),\ 0\leq x\leq \frac 1{2k},\ee with strict inequality at least at one point. Then $a=\frac 1{2k}<S$ by Proposition \ref{ex2}, and the transport (to the left) occurs. If inequality \eqref{inex1} has opposite sign, the transport is to the right. \end{example}

\begin{example} \label{ex32} Let $\nu$ be symmetric, i.e. $\nu(x)=\nu(1-x)$. Then periodicity of $\nu$ and Example \ref{ex1} imply that $S=\frac 1{2k}$.
The motor effect is provided by the condition $a\neq \frac 1{2k}$, i.e. the potential should be asymmetric, and the direction of transport is determined by the location of $a$ with respect to $\frac 1{2k}$, as for the one-state flashing ratchet \cite{vorr}. \end{example}

The rate of attachment $\nu$ of free motor proteins to the cargo is determined by the coupling between the conformation and the underlying chemical reactions. A representative situation is when the function $\nu(x)$ is peaked on narrow regions corresponding to the active sites in the conformation space, and almost vanishes outside of that regions \cite[p. 55]{jul}. Due to periodicity, in each well $[x_i,x_{i+1}]$ there is one active region. The following proposition makes very transparent the application of Theorem \ref{mainthmr} to this framework.

\begin{prop} \label{deltaf} Let $\phi_n:[A,B]\to(0,+\infty)$ be a sequence of continuous functions. Let $s_n$ be the corresponding $\phi_n$-diffusive means of $A$ and $B$. Assume that there exist numbers $s_*\in(A,B)$ and $K_*>0$ such that $\phi_n \to K_*\delta_{s_*}$ weakly-* in $C^*[A,B]$ as $n\to +\infty$. Then $s_n\to s_*$. \end{prop}

Consequently, if the active sites are located around the points $s_i=s_*+x_i,\ i=1,\dots, k,$ then the $\nu/\varsigma$-diffusive mean of $x_1$ and $x_2$ is approximately equal to $s_*$. Thus, the direction of transport is determined by the sense of the vector $\overrightarrow{s_* a}$. Hence, the factor which plays the crucial role in the transport is not the asymmetry of the potential but the location of the minima of the potential outside of the active regions. To avoid reader's confusion, we again point out that here we are referring only to the active sites where the attachment of the motors occurs.

However, the detachment of the motors may also be localized within small active zones. Normally \cite{modbio,jul}, it is assumed that the areas of detachment activity are located around the minima of the potential. Moreover, it is generally supposed that the asymmetry of the potential should correlate with the transport, i.e. the position of the minimum of the potential should be biased from the center of each well in accordance with the intended bulk movement direction, as in Example \ref{ex32}. We have just seen that all these hypotheses are redundant (provided the ratchet is well-tuned in the sense of Theorem \ref{mainthmr}).

\section {Squeezing ratchet}

We divide the proof of Theorem \ref{mainthmr} into two main steps. Firstly, in this section, we consider a somehow limiting case of the randomly flashing ratchet. It appears as a pair of purely mathematical statements, which, however, have a nice physical interpretation. Namely, we will propose an idealized ratchet-like engine, and the just mentioned mathematical facts will entail that the direction of transport provided by this machine is determined by the location of $a$ with respect to $S$.  Secondly, in the next section, we will perform some passages to the limit which will imply the claimed behaviour of the randomly flashing ratchet.

Let $G(x,y)$ be Green's function of the Sturm-Liouville operator $$\mathcal{L}= - \varsigma \frac {d^2} {dx^2}+\nu$$ on $(0,1)$ with homogeneous Neumann boundary condition.
Let $q_i(x)=G(x,a_i)$, i.e.
\be \label{drt1}\left\{\begin{array}{ll}
\nu q_i - \varsigma {q_i}_{xx}=  \delta_{a_i},
& x\in (0,1), \\
{q_i}_{x}(0)={q_i}_{x}(1)=0,
\end{array}\right.\ee and let
\be  q =  \sum\limits_{i=1}^k q_{i}.\ee As we will observe in the proof of Proposition \ref{pr1} (Appendix A), \be q_i(x)>0,\ x\in[0,1].\ee

The next theorem is our core tool for catching the motor effect.

\begin{theorem} \label{drat} There exists a unique function $Q\in C[0,1]$ solving the following problem:
\be \label{drt}\left\{\begin{array}{ll}
\nu Q - \varsigma Q_{xx}=  \sum\limits_{i=1}^k \widehat{(\nu Q)}_i \delta_{a_i},
& x\in (0,1), \\
Q_{x}(0)=Q_{x}(1)=0,\\
\int\limits_{0}^1 \nu(x) Q(x) \, dx=1.
\end{array}\right.\ee
Moreover, if \be \label{qin} q(x)\geq q(x+1/k) + \gamma, \ 0\leq x\leq 1-1/k, \ee with some $\gamma>0$, then
\be \label{qin1} Q(x)\geq Q(x+1/k) + M \gamma, \ 0\leq x\leq 1-1/k, \ee where $M=\min_{i=1,\dots,k} \widehat{(\nu q_i)}_k$.
\end{theorem}

\begin{remark} Theorem \ref{drat} may be considered as a continuous version of a purely algebraic fact, Lemma 3.2 by \citeasnoun{vorr}. \end{remark}

Our ``left drift'' condition $a<S$ arises from the following lemma.

\begin{lemma} \label{pokaz} If $a<S$, then there is $\gamma>0$ such that \eqref{qin} holds true. \end{lemma}

The results of this section can be interpreted as follows. Consider a semidiscretized device which we refer to as the ``squeezing ratchet", and which acts as follows. The particles can be in two states, ground and excited, and $\nu(x)$ and $\eta(x)$ are probabilities of transition from the first to the second state and back, resp.  The particles in the ground state diffuse with diffusion coefficient $\varsigma$. If a particle positioned at the segment $(x_i,x_{i+1})$ suddenly changes its state from the ground to the excited one, then it instantly jumps to the point $a_i$ (located to the left from the $\nu/\varsigma$-diffusive mean of the points $x_i$ and $x_{i+1}$). Then, given any initial allocation of particles, the renormalized eventual distribution $Q$ of ground particles satisfies \eqref{qin1}, i.e. their mass is transported to the left. Moreover, the excited particles are eventually concentrated at the points $a_i$, and one can observe that the asymptotic amounts $\chi_i\sim {\widehat {(\nu Q)} _i}$ of excited particles at the points $a_i$ decrease from left to right. A reflection argument shows that if $a_i$ are located to the right from the $\nu/\varsigma$-diffusive means of the corresponding endpoints $x_i$ and $x_{i+1}$, then both ground and excited mass is transported to the right.

\section{Asymptotics of the time-discretized Fokker-Planck equation and behaviour of the randomly flashing ratchet}

This section contains a series of limiting procedures which, together with the results of the previous section, will lead to derivation of Theorem \ref{mainthmr}. The reader not interested in the details of the proof of the main theorem may skip this section.

Denote by $d$ the Wasserstein metric of order two on the space of probability measures on $[0,1]$ \citeaffixed{kind3}{see e.g.}.
The convergence in Wasserstein metric is equivalent to the weak-* convergence of probability measures:\be\label{prob}d(\mu_n,\mu)\to 0 \Leftrightarrow \langle \mu_n-\mu,f\rangle\to 0,\ f\in C [0,1].\ee

Set $b(x)=\psi_x(x)/\eta(x)$. Note that $b(x)$ is zero at the extrema $a_i$ and $x_i$ of the potential $\psi$, is negative for $x_i<x < a_i$, and is positive for  $a_i<x < x_{i+1}$,  $i=1,\dots,k$. We also recall that $\kappa$ is a physical parameter inversely proportional to temperature.

Each of the following results will be used for the proof of the subsequent ones. We begin with a density-based lemma.

 \begin{lemma} \label{co0} If $\omega\in C^*[0,1]$ satisfies \be\label{om}\langle \omega, \varphi+\kappa b\varphi_x\rangle=0\ee
 for any $\varphi\in C^2[0,1]$, $\varphi_x(0)=\varphi_x(1)=0$, then $\omega=0$.  \end{lemma}

Let $\mathcal{R}$ be any uniformly integrable \cite{dies} subset of the set $$\left\{R\in L_1(0,1), \ R\geq 0,\ \int\limits_0^1 R(x)\, dx=1\right\}.$$

We now study the asymptotic behaviour of the system
\be \label{rt} \left\{\begin{array}{ll}
r -\kappa(br)_x=  R,
& x\in (0,1)\\ r\geq 0, \end{array}\right. \ee
where $R\in\mathcal{R}$ is prescribed. Clearly, any solution $r$ of \eqref{rt} satisfies $\int\limits_0^1 r(x)\, dx=1,$ so $r$ can be considered as a probability measure.

 \begin{lemma} \label{co1} There exists a unique solution to \eqref{rt}, and we have \be \lim\limits_{\kappa \to +\infty}^{}\sup\limits_{R\in \mathcal{R}}^{} d(r,  \sum\limits_{i=1}^k\widehat{R}_i\delta_{a_i}) =0.\ee  \end{lemma}

Assume that for every $\kappa>0$, $\sigma>0$ there is $R_{\sigma,\kappa}\in \mathcal{R}$ such that the problem
  \be \label{rtw}\left\{\begin{array}{ll}
w_{\sigma,\kappa} - \sigma(w_{\sigma,\kappa}/\eta)_{xx}-\kappa(bw_{\sigma,\kappa})_x=  R_{\sigma,\kappa},
& x\in (0,1), \\
(w_{\sigma,\kappa}/\eta)_{x}(0)=(w_{\sigma,\kappa}/\eta)_{x}(1)=0, \\
w_{\sigma,\kappa}\geq 0
\end{array}\right.\ee has a solution $w_{\sigma,\kappa}$. We want to study the asymptotics of the time-discretized Fokker-Planck problem \eqref{rtw} for large $\kappa$ and small $\sigma$.

 \begin{lemma} \label{co2} For each $\kappa$ there exists $\epsilon_\kappa>0$ so that \be\label{mc1} \lim\limits_{\kappa \to +\infty,\ \sigma\leq \epsilon_\kappa}^{} d(w_{\sigma,\kappa},  \sum\limits_{i=1}^k\widehat{(R_{\sigma,\kappa})}_i\delta_{a_i}) =0.\ee  \end{lemma}

We are in a position to obtain the main result of this section, namely, to prove that the randomly flashing ratchet behaves similarly to the squeezing ratchet.

  \begin{lemma} \label{co3} For each $\kappa$ there exists $\epsilon_\kappa>0$ so that the corresponding solutions of \eqref{eqb3} have the following properties: \be \label{mc2}\lim\limits_{\kappa \to +\infty,\ \sigma\leq \epsilon_\kappa}^{}\sup\limits_{0\leq x\leq 1}^{}|P(x) -Q(x)|=0,\ee \be \label{mc3} \lim\limits_{\kappa \to +\infty,\ \sigma\leq \epsilon_\kappa}^{}\widehat p _i= \frac {\widehat {(\nu Q)} _i}{\eta(a)}.\ee \end{lemma}

The proof of Theorem \ref{mainthmr} is now straightforward. Inequality \eqref{qin1} yields \be\frac {\widehat {(\nu Q)} _i}{\eta(a)}\geq \frac {\widehat {(\nu Q)} _{i+1}}{\eta(a)}+C, \ i=1,\dots,k-1.\ee Therefore, \eqref{mainr} and \eqref{qinm} are direct consequences of Theorem \ref{drat} and of Lemmas \ref{pokaz} and \ref{co3}.

  \section{Discussion} We investigate the randomly flashing ratchet model for motor proteins. In this model, every Brownian particle can be in two states: with or without a motor protein. The particles carrying a motor protein are sensitive to a chemically-induced periodic potential. The motor-free particles diffuse normally. The switch between the two states happens when a particle loses or seizes a motor protein, and the probabilities of these events are prescribed.

The real nature of transport generated by motor proteins is much more complex, but the model under consideration is flexible enough to encompass the four main factors involved: chemical potential, conformation, attachment/detachment of motor proteins, and diffusion.

We have rigorously shown that unidirectional transport of mass occurs when the diffusion of the potential-sensitive particles is strongly dominated by the influence of the potential. The direction of transport is unambiguously determined by a certain interrelation between the asymmetries of the potential and of the function describing the probability of attachment of motor protein molecules to the cargo (in other words, of the rate of transition from the potential-insensitive to the potential-sensitive state): in particular, at least one of them should be asymmetric to secure the transport effect.

In the literature it is often assumed \cite{jul} that the attachment of motor proteins to transportable objects predominantly happens within a narrow active region determined by the conformation. This region is unique in every period of the potential. As an application of our general result, we show that in this specific framework the ratchet mechanism moves the transported matter to the left/right provided  the minimum of the potential  within each period is located to the left/right of the attachment zone. We do not require additional hypotheses imposed in earlier biophysical papers such that correlation of the asymmetry of the potential with the intended direction of transport or localization of the detachment areas around the minima of the potential.

Our theory is also able to detect the transport direction in the case when there are two or more attachment sites within each period or their spatial distribution becomes more homogeneous.

Another popular model for molecular motors is the classical or deterministic flashing ratchet. In the biomotor context, it may be interpreted assuming that during some period of time all the particles are bound to motor proteins, and are therefore potential-sensitive. Then, suddenly, all the tiny ``parcels'' simultaneously lose their carriers and start to diffuse freely. After a certain time interval, the motors synchronously reattach to their cargoes, and the whole process repeats on and on. Thus, conformation is completely disregarded. The direction of transport is dictated by the location of the minima of the asymmetric periodic chemical potential with respect to the centres of the potential wells \cite{vorr}. The deterministic flashing ratchet has promising applications in many areas of science. However, in connection with motor proteins, the random ratchet seems to be more pliable and realistic.

Our analysis shows that predictions of the random and deterministic flashing ratchets may significantly differ. For instance, if the region of attachment activity is unique (up to spatial periodicity), and the minimum of the chemical potential is located between that active zone and the centre of the potential well, then the bulk transport directions anticipated by the two models are exactly the opposite of one another. On the other hand, when the distribution of attachment sites is symmetric (no matter how many active areas are present in each well), the forecasts of the models essentially coincide. These observations give a new insight for experimental verification and validation of the flashing ratchets.

We now summarize the major biological implications of the paper. Our theoretical analysis has highlighted the main factors that govern the operation of biomotors: chemical potential, architecture of the involved proteins, structural changes in the system, including  attachment/detachment of motor proteins, and diffusion. We have illustrated that the interplay between these factors is very complex, and their subtle alterations may change behaviour of the system, in particular, the bulk transport direction.  Our key finding is that the spatial distribution of the active sites (where the attachment of motor proteins occurs) is as important as the form of the potential. Previous rigorous studies revolved around the primacy of asymmetry of the chemical potential, whereas the influence of distribution of attachment/detachment sites was considered as a subordinate and tributary factor or disregarded at all. We have also made it possible to predict the collective behaviour of motor proteins and their cargoes in vivo and in vitro. We have shed more light on the understanding of the mechanism of motor proteins, and have showed that the value of qualitative mathematical analysis cannot be underrated prior to designing biomotor-powered hybrid devices and artificial biomotors.
\appendix

\section{Proof of Propositions \ref{pr1} and \ref{ex2}}

Let $G(x,y)$ be Green's function of the Sturm-Liouville operator $$\mathcal{L}= - \frac {d^2} {dx^2}+\phi$$ on $(A,B)$ with homogeneous Neumann boundary condition. Then $U(x)=G(x,s)$ is a solution to \eqref{defm} if an only if $G(A,s)=G(B,s)$.

By the distributional maximum principle \cite[Theorem B]{litt}, $G(x,y)>0$. Observe that \be G_x(x,A)=\int\limits_{B}^x\phi(z) G(A,z) \, dz,  \ee
\be G_x(x,B)=\int\limits_{A}^x\phi(z) G(B,z) \, dz.  \ee
 Hence, the function $G(x,A)$ is decreasing in $x$, and $G(x,B)$ is increasing.
 Thus, the function $$g(x)=G(A,x)-G(B,x)$$ is also (strictly) decreasing. At the ends of the segment, we have $g(A)=G(A,A)-G(B,A)> 0$ and $g(B)=G(A,B)-G(B,B) < 0$. Since $g$ is a continuous function, there is unique $s\in (A,B)$ such that $G(A,s)=G(B,s)$. The first proposition is proven.

Now,  let $$g_1(x)=G(\frac {A+B} 2,x)-G(\frac {A+B} 2,A+B-x).$$ We claim that \be g_1(A)>0.\ee

Since $G$ is Green's function, \be (g_1)_x(A)=0,\ee and \begin{multline} (g_1)_{xx}(x)\\ =\phi(x)G(\frac {A+B} 2,x)-\phi(A+B-x)G(\frac {A+B} 2,A+B-x)\\ =\phi(x) g_1(x)+[\phi(x)-\tilde\phi(x)]G(\frac {A+B} 2,A+B-x)\\ \leq \phi(x) g_1(x), \ A\leq x< \frac {A+B} 2,\end{multline} and the inequality is strict at least at one point. In particular, $g_1$ cannot be identically zero.

 Assume that $g_1(A)\leq 0$. By the maximum principle, $g_1$ cannot have non-positive minima within $(A,\frac {A+B} 2)$. But \be\label{go} g_1(\frac {A+B} 2)=0,\ee so $A$ must be a minimum point.  Let $g_2(x)=g_1(x)-g_1(A)$. Then $g_2$ is non-negative, and \be (g_2)_{xx}(x) \leq \phi(x) g_2(x), \ A\leq x< \frac {A+B} 2.\ee Thus, \be (g_2)_{x}(x) \leq \int\limits_A^x \phi(t) g_2(t)\,dt, \ A\leq x\leq \frac {A+B} 2.\ee  Since $g_2(A)=0$, by the mean value theorem, \begin{multline}  g_2(x)=g_2(x)-g_2(A)= (x-A)(g_2)_x(c) \\ \leq (x-A)\int\limits_A^c \phi(t) g_2(t)\,dt \leq \frac 12 (B-A)\int\limits_A^x \phi(t) g_2(t)\,dt\end{multline} for some $c$, $A<c<x\leq \frac {A+B} 2$. The Gronwall lemma implies $g_2\equiv 0$, so, by \eqref{go}, $g_1\equiv 0$, and we get a contradiction.

Hence, $$g(\frac {A+B} 2)=g_1(A)> 0.$$ We have observed above that $g(B)<0$, so there is $s\in (\frac {A+B} 2,B)$ such that $g(s)=0$, and this number $s$ is the $\phi$-diffusive mean.

\section{Proof of Proposition \ref{deltaf}}

Without loss of generality, $[A,B]=[0,1]$.
Denoting the solutions of \eqref{defm} corresponding to $\phi_n$ and $s_n$ by $U_n$, we infer that \be \label{de1} \int\limits_0^1 \phi_n(x) U^2_n(x) +(U_n)_x^2\, dx = U_n(s_n).\ee The function $U_n$ is positive and thus convex on $(0,s_n)$ and $(s_n,1)$. Hence, due to the boundary conditions, the derivative $(U_n)_x(x)$ is positive/negative when $x< s_n$ / $x >s_n$, and tends toward its supremum/infimum as $x$ approaches $s_n$ from the left/right, resp.  Moreover, the maximum of $U_n$ is achieved at $s_n$, and the minima are reached at $0$ and $1$.

It is easy to see that \be\min\limits_{x\in[0,1]} U_n^2(x)\leq \frac 2 {K_*} \int\limits_0^1 \phi_n(x) U^2_n(x)\,dx\ee for large $n$, and \be \max\limits_{x\in[0,1]} U_n^2(x)\leq 2\left[\min\limits_{x\in[0,1]} U_n^2(x)+\int\limits_0^1 (U_n)_x^2\, dx\right].\ee Hence, \eqref{de1} implies
\be \|U_n^2\|_{C[0,1]} \leq C \|U_n\|_{C[0,1]},\ee whence \be \label{unest} \|U_n\|_{C[0,1]} \leq C.\ee

Without loss of generality, there exists a limit $s^*\in[0,1]$ of the sequence $\{s_n\}$. Assume that $s^*\neq s_*$, or, more particularly, $s_*>s^*$ (the opposite case may be examined in a similar way). This yields that the value of the integral $\int\limits_0^{s_n} \phi_n(x) $ tends to zero as $n\to + \infty$. But integration of \eqref{defm} implies that the left derivative $(U_n)^\prime_- (s_n)$ is equal to $\int\limits_0^{s_n} \phi_n(x) U_n(x)\,dx$. Due to \eqref{unest}, this integral goes to zero as $n\to + \infty$. Using the information on the behaviour of $U_n$ summarized after equality \eqref{de1}, we conclude that $\|(U_n)_x\|_{C[0,s_n]}\to 0$. Due to \eqref{unest}, without loss of generality there exists a constant $U_*$ such that \be \max\limits_{x\in[0,1]} U_n(x) \to U_*.\ee But \begin{multline}\max\limits_{x\in[0,1]} U_n(x)- \min\limits_{x\in[0,1]} U_n(x)=U_n(s_n)-U_n(0)\\ \leq s_n\|(U_n)_x\|_{C[0,s_n]}\to 0.\end{multline} Hence, \be U_n(x) \to U_*\ee uniformly on $[0,1]$.

Let us now test \eqref{defm} with a smooth function $h$ such that $h(s_*)=0$, $h(s^*)=1$, $h_x(0)=0$, $h_x(1)=0$, obtaining \be \int\limits_0^1 \phi_n(x) U_n(x) h(x)-U_n(x) h_{xx}(x)\, dx = h(s_n).\ee Passing to the limit, we find that the left-hand side goes to zero and the right-hand side goes to one, arriving at a contradiction.

\section{Proof of Theorem \ref{drat}}

Consider the set of functions \be \notag B=\left\{y(x)\in L_1(0,1)\,\middle|\, \int\limits_{0}^1 \nu(x) y(x) \, dx=1, \ \widehat{(\nu y)}_k \geq M, \right.\ee \be \notag
 \left.\phantom{\int\limits_{i=1}^m} y(x)\geq y(x+1/k)+M \gamma , \ \mathrm{for\ a.a.}\ 0\leq x\leq 1-1/k.\right\} \ee

Inverse induction shows that for any $i=1,\dots, k$ and $y\in B$ one has \be \label{iin} \widehat{(\nu y)}_i \geq M. \ee

Let us define a mapping $\mathcal{A}$ on $B$. For each $y\in B$, we let $\mathcal{A}(y)=Y$, where $Y$ is the solution of the problem
\be \label{drt2}\left\{\begin{array}{ll}
\nu Y - \varsigma Y_{xx}=  \sum\limits_{i=1}^k \widehat{(\nu y)}_i \delta_{a_i},
& x\in (0,1), \\
Y_{x}(0)=Y_{x}(1)=0.
\end{array}\right.\ee
To put it differently,
\be \label{drt3} Y=\sum\limits_{i=1}^k \widehat{(\nu y)}_i q_{i}.\ee

Then, the set $B$ is invariant for the map $\mathcal{A}$. In fact, let $y\in B$. Then \eqref{drt2} implies
$$\int\limits_{0}^1 \nu(x) Y(x) \, dx= \sum\limits_{i=1}^k \widehat{(\nu y)}_i = \int\limits_{0}^1 \nu(x) y(x) \, dx=1.$$
Further, $$\widehat{(\nu Y)}_k=\sum\limits_{i=1}^k \widehat{(\nu y)}_i \widehat{(\nu q_i)}_k\geq M\sum\limits_{i=1}^k \widehat{(\nu y)}_i =M.$$ Finally, fix $x_*\in[0,1-1/k]$. Then there is a number $n$ such that $x_*\in[x_n,x_{n+1})$. Set \be N_*=\left\{\begin{array}{ll} \widehat{(\nu y)}_n,
& x_*\leq a_n,\\
\widehat{(\nu y)}_{n+1}, & x_* > a_n.
\end{array}\right.\ee We claim that \be\label{ny} (\widehat{(\nu y)}_i - N_*)[q_i(x_*)-q_i(x_*+1/k)]\geq 0,\ i=1,\dots,k. \ee
Indeed, integration of \eqref{drt1} gives \be\varsigma {q_i}_{x}(x)=\int\limits_{0}^x\nu(z) q_i(z) \, dz, \ x< a_i,  \ee
\be\varsigma {q_i}_{x}(x)=\int\limits_{1}^x\nu(z) q_i(z) \, dz, \ x> a_i.  \ee
Thus, the function $q_i$ is increasing on the segment $[0,a_i]$ and decreasing on $[a_i,1]$. Assume first $x_*\leq a_n$.
Then, if $i<n$, we have $q_i(x_*)>q_i(x_*+1/k)$.  Since $y\in B$, we also have  $\widehat{(\nu y)}_i >\widehat{(\nu y)}_n= N_*$, and \eqref{ny} holds true. If $i>n$, we have $q_i(x_*)<q_i(x_*+1/k)$ and $\widehat{(\nu y)}_i < N_*$, and \eqref{ny} again holds. If $i=n$, \eqref{ny} is trivial. Now, let $x_*> a_n$. In this case, if $i<n+1$,  $q_i(x_*)>q_i(x_*+1/k)$ and $\widehat{(\nu y)}_i >\widehat{(\nu y)}_{n+1}= N_*$; if $i>n+1$, $q_i(x_*)<q_i(x_*+1/k)$ and $\widehat{(\nu y)}_i < N_*$; and if $i=n+1$, \eqref{ny} is again trivial.

Formulas \eqref{drt3}, \eqref{ny}, \eqref{qin} and \eqref{iin} yield $$Y(x_*)- Y(x_*+1/k) = \sum\limits_{i=1}^k \widehat{(\nu y)}_i [q_i(x_*)-q_i(x_*+1/k)]$$ $$\geq \sum\limits_{i=1}^k N_* [q_i(x_*)-q_i(x_*+1/k)]= N_* [q(x_*)-q(x_*+1/k)]\geq M\gamma,$$
so the invariance of $B$ is confirmed.

Observe that $\mathcal{A}$ is a compact linear operator in $L_1(0,1)$. Indeed, let $\mathbb{B}$ be the unit ball of the space $L_1(0,1)$. Due to \eqref{drt3}, its image $\mathcal{A}(\mathbb{B})$ is a bounded subset of the linear span of $\{q_1,\dots,q_k \}$, thus being a relatively compact subset of a finite-dimensional subspace of $L_1(0,1)$.

Let us show that \eqref{drt} may have at most one solution, so $\mathcal{A}$ can have at most one fixed point in $B$. If not, let ${\tilde{Q}}$ be the difference of two distinct solutions. Then \be\label{sumc}\int\limits_{0}^1 \nu(z) \tilde{Q}(z) \, dz=0.\ee Moreover, \be \label{drt4} {\tilde{Q}}=\sum\limits_{i=1}^k \widehat{(\nu {\tilde{Q}})}_i q_{i},\ee whence \be \label{eq22} \widehat{(\nu {\tilde{Q}})}_j=\sum\limits_{i=1}^k \widehat{(\nu {\tilde{Q}})}_i \widehat{(\nu q_i)}_j,\ j=1,\dots,k.\ee From \eqref{drt1} we deduce \be\int\limits_{0}^1 \nu(z) q_i(z) \, dz=1.\ee Therefore, the matrix $[\mathcal{P}_{ij}]=\left[\widehat{(\nu q_i)}_j \right]$ is \emph{ergodic}, i.e. it has positive entries, and the sum of the elements in every row is equal to one. By the Perron-Frobenius theorem, it has an eigenvector $[\xi_i]$ corresponding to the simple eigenvalue $1$, so that $\xi_j=\sum\limits_{i=1}^k \xi_i \mathcal{P}_{ij}$, and all the components $\xi_i$ are positive. On the other hand, by \eqref{eq22}, $[\Xi_i]=\left[\widehat{(\nu \tilde{Q})}_i \right] $ is another eigenvector of $[\mathcal{P}_{ij}]$ corresponding to the same eigenvalue. The sum of its components is zero due to \eqref{sumc}, so it cannot be collinear with $[\xi_i]$ unless it is a zero vector. Since $1$ is a simple eigenvalue, all $\widehat{(\nu \tilde{Q})}_i$ are zeros, so $\tilde{Q}\equiv 0$ by virtue of \eqref{drt4}.

The set $B$ is closed, convex and bounded in $L_1(0,1)$. By Schauder's fixed point principle, $\mathcal{A}$ has a fixed point $Q$ in $B$, which is automatically a solution to \eqref{drt}.  It remains to notice that $Q$ is continuous as a linear combination of $q_i$, so \eqref{qin1} holds for all $0\leq x\leq 1-1/k$.

\section{Proof of Lemma \ref{pokaz}}

Let us notice that \be q=u+v,\ee where $u$ and $v$ are the (unique) solutions to the following problems
\be \label{r1}\left\{\begin{array}{ll}
\nu u - \varsigma {u}_{xx}=  \sum\limits_{i=1}^k\delta_{S_i},
& x\in (0,1), \\
{u}_{x}(0)={u}_{x}(1)=0,
\end{array}\right.\ee
\be \label{r2}\left\{\begin{array}{ll}
\nu v - \varsigma {v}_{xx}=  \sum\limits_{i=1}^k(\delta_{a_i}-\delta_{S_i}),
& x\in (0,1), \\
{v}_{x}(0)={v}_{x}(1)=0,
\end{array}\right.\ee and \be S_i= S+x_i,\ i=1,\dots, k.\ee
Since $S$ is the $\nu/\varsigma$-diffusive mean of $0$ and $1/k$, there exists a solution $U_1$ to the problem
\be \label{r3}\left\{\begin{array}{ll}
\nu U_1 - \varsigma {U_1}_{xx}=  \delta_{S_1},
& x\in (0,1/k), \\
{U_1}_{x}(0)={U_1}_{x}(1/k)=0,\\
U_1(0)=U_1(1/k).
\end{array}\right.\ee
The solution $u$ to \eqref{r1} can be constructed in the following way: \be u(x)=U_1(x-x_i),\ x_i\leq x\leq x_i+1/k,\ i=1,\dots, k.\ee Thus, $u$ is $1/k$-periodic, i.e. \be \label{qinu} u(x)- u(x+1/k)=0, \ 0\leq x\leq x_k. \ee
Set \be \label{qinu1} d(x)=v(x)- v(x+1/k), \ 0\leq x\leq x_k. \ee Then it suffices to show that \be\gamma=\min_{0\leq x\leq x_k} d(x)>0.\ee
Note that
\be \label{r4}
\nu d - \varsigma {d}_{xx}=  0,
\ x\in (0,x_k).\ee By the maximum principle, if the minimum of $d$ is non-positive, it is attained at $0$ or $x_k$. To ascertain that this cannot happen, we are going to prove that  \be d_x(0)<0, \ d_x(x_k)>0.\ee

Set \be V(x)=\int\limits_{0}^x\nu(z) v(z) \, dz,\ee and let $\theta$ be the solution of the Cauchy problem
\be \label{theta}\left\{\begin{array}{ll}
\theta_x=  \sum\limits_{i=1}^k(\delta_{a_i}-\delta_{S_i}),
& x\in (0,1), \\
\theta(0)=0. \end{array}\right.\ee
Note that $\theta$ is non-negative and $1/k$-periodic.

Integration of \eqref{r2} gives
\be \label{f40} V-\varsigma v_x= \theta.\ee
Therefore
\be d_x=\frac {V(x)-V(x+1/k)}{\varsigma}.\ee
From \eqref{f40} we deduce
\be \label{f41} V(0)=V(1)=0,\ee and
\be \label{f43} V-\varsigma \left(\frac {V_x}{\nu}\right)_x= \theta\geq 0.\ee Using the distributional maximum principle \cite[Theorem B]{litt}, we conclude that \be \label{f42} V(x)> 0, \ 0<x<1,\ee so
 \be d_x(0)=-\frac {V(1/k)}{\varsigma}<0, \ d_x(x_k)=\frac {V(1-1/k)}{\varsigma}>0.\ee

\section{Proof of Lemma \ref{co0}} It suffices to prove that the set $$O=\left\{\varphi+\kappa b\varphi_x\Big| \varphi\in C^2[0,1], \varphi_x(0)=\varphi_x(1)=0\right\}$$ is dense in $C[0,1]$.

Let $h\in C^2[0,1]$ be an arbitrary function which is locally constant near the zeros of $b$. These functions constitute a dense subset $O_1$ of $C[0,1]$.  Let \begin{multline}\varphi(x)=h(x)+ \int\limits_{x}^{a_i} \exp\left(\int\limits_x^y\frac 1 {\kappa b(t)}\,dt \right)h_y(y)\,dy,\\ x_i< x< x_{i+1}, \ i=1,\dots, k, \\ \varphi(x_i)=h(x_i),\ i=1,\dots, k+1.\end{multline} Clearly, $\varphi$ is equal to a constant $c^-_i$ (resp. $c^+_i$) in a left (resp. right) neighbourhood of the point $x_i$. But \be \label{f65} h=\varphi+\kappa b\varphi_x,\ee so $c^-_i=c^+_i=h(x_i)$. Thus, $\varphi$ is $C^2$-smooth and $\varphi_x(0)=\varphi_x(1)=0$. By virtue of \eqref{f65}, $O_1$ is contained in $O$.

\section{Proof of Lemma \ref{co1}}

The solution $r$ to \eqref{rt} can be written explicitly: \be r(x)=-\frac 1 {\kappa b(x)}\int\limits_{x_i}^x \exp\left(\int\limits_s^x\frac 1 {\kappa b(t)}\,dt \right)R(s)\,ds,\ x_i<x < a_i,\ee
 \be r(x)=\frac 1 {\kappa b(x)}\int\limits_{x}^{x_{i+1}} \exp\left(\int\limits_s^x\frac 1 {\kappa b(t)}\,dt \right)R(s)\,ds,\ a_i<x < x_{i+1}.\ee  Note that it is unique in $L_1(0,1)$. Indeed, if $r_1$ is another solution, then $\omega=r-r_1$ satisfies the conditions of Lemma \ref{co0}.

 We need to show that $r\to  \sum\limits_{i=1}^k\widehat{R}_i\delta_{a_i}$ weakly-$*$, uniformly with respect to $R\in\mathcal{R}$. It suffices to prove that, for each $i$, $r\to \widehat{R}_i\delta_{a_i}$ weakly-$*$ on the interval $(x_i,x_{i+1})$, uniformly in $R\in\mathcal{R}$. We restrict ourselves to the case $i=1$, and the others are analogous.

 We calculate, integrating by parts,  \begin{multline} \widehat{r}_1=\int\limits_0^{1/k} r(x)\,dx \\=-\int\limits_0^{a} \frac 1 {\kappa b(x)}\int\limits_{0}^x \exp\left(\int\limits_s^x\frac 1 {\kappa b(t)}\,dt \right)R(s)\,ds\,dx \\+\int\limits_{a}^{1/k}\frac 1 {\kappa b(x)}\int\limits_{x}^{1/k} \exp\left(\int\limits_s^x\frac 1 {\kappa b(t)}\,dt \right)R(s)\,ds\,dx\\
 =\left[\int\limits_{0}^x \exp\left(\int\limits_s^x\frac 1 {\kappa b(t)}\,dt \right)R(s)\,ds\right]_a^0+\int\limits_0^{a} R(x)\,dx\\ +\left[\int\limits_{1/k}^x \exp\left(\int\limits_s^x\frac 1 {\kappa b(t)}\,dt \right)R(s)\,ds\right]_{1/k}^a+\int\limits_a^{1/k} R(x)\,dx\\= \widehat{R}_1-\int\limits_{0}^{1/k} \exp\left(\int\limits_s^a\frac 1 {\kappa b(t)}\,dt \right)R(s)\,ds= \widehat{R}_1.\end{multline}

 Let us show that for every $x_*\in(0,a)$ \be\label{ll}\lim\limits_{\kappa \to +\infty}^{}\int\limits_0^{x_*} r(x)\,dx= 0,\ee
 uniformly in $R\in\mathcal{R}$. Indeed, let $s_\kappa<x_*$ be such that \be\int\limits_{x_*}^{s_\kappa}\frac 1 { b(t)}\,dt= \sqrt\kappa.\ee Observe that $s_\kappa\to 0$ as $\kappa\to +\infty$. We have \begin{multline}\int\limits_0^{x_*} r(x)\,dx\\=-\int\limits_0^{x_*} \frac 1 {\kappa b(x)}\int\limits_{0}^x \exp\left(\int\limits_s^x\frac 1 {\kappa b(t)}\,dt \right)R(s)\,ds\,dx\\
 =\left[\int\limits_{0}^x \exp\left(\int\limits_s^x\frac 1 {\kappa b(t)}\,dt \right)R(s)\,ds\right]_{x_*}^0+\int\limits_0^{x_*} R(x)\,dx\\
 =\int\limits_{0}^{x_*} \left[1-\exp\left(\int\limits_s^{x_*}\frac 1 {\kappa b(t)}\,dt \right)\right]R(s)\,ds\\ \leq \int\limits_{0}^{s_\kappa} \left[1-\exp\left(\int\limits_s^{x_*}\frac 1 {\kappa b(t)}\,dt \right)\right]R(s)\,ds + \int\limits_{s_\kappa}^{x_*} \left[1-\exp\left(\int\limits_{s_\kappa}^{x_*}\frac 1 {\kappa b(t)}\,dt \right)\right]R(s)\,ds\\ \leq
 \int\limits_{0}^{s_\kappa} R(s)\,ds +  [1-\exp(-\kappa^{-1/2})]\int\limits_{s_\kappa}^{x_*}R(s)\,ds\to 0 \end{multline} as $\kappa\to +\infty$.

Due to \eqref{ll}, for any $f_0\in C[0,1/k]$,  \be\label{ll1}\lim\limits_{\kappa \to +\infty}^{}\int\limits_0^{x_*} r(x)f_0(x)\,dx= 0,\ee
 uniformly in $R\in\mathcal{R}$.

 Similarly, for all $x^*\in(a,1/k)$ and $f_0\in C[0,1/k]$,  \be\label{ll2}\lim\limits_{\kappa \to +\infty}^{}\int\limits_{x^*}^{1/k} r(x)f_0(x)v= 0,\ee
 uniformly in $R\in\mathcal{R}$.

 Fix $\varepsilon>0$ and $f\in C[0,1/k]$. Let $x_*$ and $x^*$ be so close to $a$ that $|f(x)-f(a)|\leq \varepsilon/2$ provided $x_*\leq x \leq x^*$. Then \be \left|\int\limits_{x_*}^{x^*}  r(x)[f(x)-f(a)]\,dx\right|\leq \varepsilon/2.\ee Due to \eqref{ll1} and \eqref{ll2} with $f_0=f-f(a)$, \be \left|\int\limits_{0}^{1/k}  r(x)[f(x)-f(a)]\,dx-\int\limits_{x_*}^{x^*}  r(x)[f(x)-f(a)]\,dx\right|\leq \varepsilon/2\ee for sufficiently large $\kappa$. Thus, \be\langle r- \widehat{R}_1\delta_{a_1},f\rangle=\langle r- \widehat{r}_1\delta_{a},f\rangle=\int\limits_{0}^{1/k}  r(x)[f(x)-f(a)]\,dx\leq\varepsilon.\ee

\section{Proof of Lemma \ref{co2}}

Let $r_{\sigma,\kappa}$ be the solution of the system
\be \label{rtk} \left\{\begin{array}{ll}
r_{\sigma,\kappa} -\kappa(br_{\sigma,\kappa})_x=  R_{\sigma,\kappa},
& x\in (0,1)\\ r_{\sigma,\kappa}\geq 0. \end{array}\right. \ee Then, by Lemma \ref{co1}, \be\lim\limits_{\kappa \to +\infty}^{} d(r_{\sigma,\kappa},  \sum\limits_{i=1}^k\widehat{(R_{\sigma,\kappa})}_i\delta_{a_i}) =0,\ee uniformly in $\sigma$.
Thus, it suffices to prove that for every $\kappa$ there is $\epsilon_\kappa>0$ such that \be\notag \lim\limits_{\kappa \to +\infty,\ \sigma\leq \epsilon_\kappa}^{} d(w_{\sigma,\kappa},  r_{\sigma,\kappa}) =0.\ee This would follow from the claim that for every $\kappa$ there is $\epsilon_\kappa>0$ so that for $\sigma\leq \epsilon_\kappa$ we have $d(w_{\sigma,\kappa},  r_{\sigma,\kappa}) < 1/\kappa$. If it is not true, then for some $\kappa$ there exists a sequence $\sigma_n\to 0$ such that $$ d(w_{\sigma_n,\kappa},  r_{\sigma_n,\kappa}) \geq 1/\kappa.$$

Since $w_{\sigma_n,\kappa}$ and $r_{\sigma_n,\kappa}$ are solutions of the problems \eqref{rtw} and \eqref{rtk}, we have \be\label{ts}- \sigma_n \langle w_{\sigma_n,\kappa},\varphi_{xx}/\eta	\rangle+\langle w_{\sigma_n,\kappa},\varphi+\kappa b\varphi_x\rangle= \langle R_{\sigma_n,\kappa},\varphi\rangle,\ee
\be\label{tbs}\langle r_{\sigma_n,\kappa},\varphi+\kappa b\varphi_x\rangle= \langle R_{\sigma_n,\kappa},\varphi\rangle,\ee for any $\varphi\in C^2[0,1]$, $\varphi_x(0)=\varphi_x(1)=0$. Since the sequences $w_{\sigma_n,\kappa}$ and $r_{\sigma_n,\kappa}$ lie in the space of probability measures, which is weakly-* compact, without loss of generality there exist their weak-* limits $w_\kappa$ and $r_\kappa$. Clearly, \be d(w_{\kappa},  r_{\kappa}) \geq 1/\kappa.\ee On the other hand, taking the difference of \eqref{ts} and \eqref{tbs}, and passing to the limit, we find $\langle w_{\kappa}-r_{\kappa},\varphi+\kappa b\varphi_x\rangle=0,$ so $w_{\kappa}=r_{\kappa}$ by Lemma \ref{co0}, and we arrive at a contradiction.

\section{Proof of Lemma \ref{co3}}

   Multiplying the second equation in \eqref{eqb3} by $P$ and integrating, we find \be\int\limits_0^1\nu(x) P^2(x) -\varsigma P_{xx}(x)P(x)\,dx=\int\limits_0^1\eta(x) p(x)P(x)\,dx,\ee
   whence \be\inf\limits_{0\leq x\leq 1}^{}\nu(x)\int\limits_0^1 P^2(x)\,dx +\varsigma \int\limits_0^1 P_{x}^2(x)\,dx\leq \sup\limits_{0\leq x\leq 1}^{}P(x).\ee
Hence,  \be\label{joe11}\|P\|^2_{W_2^1(0,1)} \leq C\|P\|_{C[0,1]}\leq C\|P\|_{W_2^1(0,1)}\leq C.\ee

The pair $(w_{\sigma,\kappa}, R_{\sigma,\kappa})=(\eta p,\nu P)$ satisfies \eqref{rtw}. Due to \eqref{joe11}, the set $\mathcal{R}=\{\nu P\,|\, \sigma>0,\kappa>0\}$ is uniformly bounded and thus uniformly integrable. By Lemma \ref{co2}, for every $\kappa$ there exists $\epsilon_\kappa>0$ such that \be\label{mc4} \lim\limits_{\kappa \to +\infty,\ \sigma\leq \epsilon_\kappa}^{} d(\eta p,  \sum\limits_{i=1}^k\widehat{(\nu P)}_i\delta_{a_i}) =0.\ee

Assume that \eqref{mc2} is not true, i.e. there exist $\delta>0$ and sequences $\kappa_n\to \infty$ and $\sigma_n\leq \epsilon_{\kappa_n}$ such that for the corresponding solutions $(p_n, P_n)=(p_{\sigma_n,\kappa_n},P_{\sigma_n,\kappa_n})$ to \eqref{eqb3} we have $\|P_{n}-Q\|_{C[0,1]}>\delta$. Since the embedding $W_2^1(0,1)\subset C[0,1]$ is compact, without loss of generality we may assume that $P_n$ converges to some limit $P_0$ in $C[0,1]$. Obviously, \be\label{pq}\|P_{0}-Q\|_{C[0,1]}\geq\delta.\ee
Passing to the limit in the second, forth and the last equations in \eqref{eqb3} --- the combination of the first two is understood in the weak sense \eqref{ewu} --- and remembering \eqref{mc4}, we find $$\nu P_0 - \varsigma (P_0)_{xx}=  \sum\limits_{i=1}^k \widehat{(\nu P_0)}_i \delta_{a_i},$$ $$ (P_0)_x(0)=(P_0)_x(1)=0,$$
$$\int\limits_0^1 \nu(x)P_0(x)\, dx =1.$$ By Theorem \ref{drat}, $P_0$ coincides with $Q$, which  contradicts \eqref{pq}.

From \eqref{mc4} we deduce \be\label{mc5} \lim\limits_{\kappa \to +\infty,\ \sigma\leq \epsilon_\kappa}^{} d\left(p,  \sum\limits_{i=1}^k\frac{\widehat{(\nu P)}_i\delta_{a_i}}{\eta(a_i)} \right)=0.\ee Due to \eqref{mc2} and $1/k$-periodicity of $\eta$, \eqref{mc5} implies that \be\label{mc6} p \to \sum\limits_{i=1}^k\frac{\widehat{(\nu Q)}_i\delta_{a_i}}{\eta(a)}.\ee
weakly-* as $\kappa \to +\infty,\ \sigma\leq \epsilon_\kappa$. Taking test functions which are equal to $1$ in one of the wells and are zero at the minima of the potential located outside of that well, we derive \eqref{mc3} from \eqref{mc6}.


\begin{thebibliography}{99}

\harvarditem{Ait-Haddou and Herzog}{2003}{modbio} Ait-Haddou R, Herzog W (2003) Brownian ratchet models of molecular motors. Cell Biochem Biophys 38:191-214.

\harvarditem{Amengual et al.}{2004}{ameng} Amengual P, Allison A, Toral R, Abbott D (2004) Discrete-time ratchets, the Fokker-Planck equation and Parrondo's paradox. Proc R Soc Lond Ser A Math Phys Eng Sci 460:2269-2284.

\harvarditem{Astumian}{1997}{ast} Astumian RD (1997) Thermodynamics and kinetics of a Brownian motor. Science 276:917-922.




\harvarditem{Chipot et al.}{2004}{chip} Chipot M, Hastings S, Kinderlehrer D (2004) Transport in a molecular motor system. M2AN Math Model Numer Anal 38:1011-1034.

\harvarditem{Chipot et al.}{2003}{chipm} Chipot M, Kinderlehrer D, Kowalczyk M (2003) A variational principle for molecular motors. Meccanica 38: 505-518.

\harvarditem{Diestel}{1991}{dies} Diestel J (1991) Uniform integrability: an introduction.
School on Measure Theory and Real Analysis (Grado, 1991).
Rend Istit Mat Univ Trieste 23:41-80.

\harvarditem{Friedman and Hu}{2007}{af1} Friedman A, Hu B (2007)
Uniform convergence for approximate traveling waves in linear reaction-hyperbolic systems.
Indiana Univ Math J 56:2133-2158.


\harvarditem{Harmer et al.}{2000}{har} Harmer G, Abbott D, Taylor P (2000)
The paradox of Parrondo's games.
Proc R Soc Lond Ser A Math Phys Eng Sci 456:247-259.

\harvarditem{Hastings et al.}{2007/08}{hast} Hastings S, Kinderlehrer D, McLeod JB (2007/08) Transport in multiple state systems. SIAM J Math Anal 39:1208-1230.

\harvarditem{Jordan et al.}{1998}{kind3} Jordan R, Kinderlehrer D, Otto F (1998) The variational formation of the Fokker-Planck equation. SIAM J Math Anal 29:1-17.

\harvarditem{J\"{u}licher}{1999}{jul} J\"{u}licher F (1999) Force and motion generation of molecular motors: A generic description. In M\"{u}ller SC, Parisi J, Zimmermann W (eds) Transport and structure, Their Competitive Roles in Biophysics and Chemistry, Lecture Notes in Physics, V 532, Springer, New York, pp 46-74.

\harvarditem{J\"{u}licher et al.}{1997}{mod2} J\"{u}licher F, Ajdari A, Prost J (1997) Modeling molecular motors. Rev Modern Phys 69:1269-1281.

\harvarditem{Kinderlehrer and Kowalczyk}{2002}{kind1} Kinderlehrer D, Kowalczyk M (2002) Diffusion-mediated transport and the flashing ratchet.
Arch Rat Mech Anal 161:149-179.

\harvarditem{Littman}{1963}{litt} Littman W (1963)
Generalized subharmonic functions: Monotonic approximations and an improved maximum principle.
Ann Scuola Norm Sup Pisa (3) 17:207-222.

\harvarditem{Mennerat-Robilliard}{1999}{cold} Mennerat-Robilliard C (1999) Atomes froids dans des r\'{e}seaux optiques - Quelques
facettes surprenantes d'un syst\`eme mod\`ele. Th\`ese de Doctorat de l'Universit\'e Paris-VI, LKB/ENS - Universit\'e Paris-VI.

\harvarditem{Mirrahimi and Souganidis}{2013}{mer} Mirrahimi S, Souganidis PE (2013) A homogenization approach for the motion of motor proteins. Nonlinear Differ Equ Appl 20:129-147.

\harvarditem{Parrondo et al.}{1998}{modmain1} Parrondo JMR, Blanco JM, Cao FJ, Brito R (1998) Efficiency of Brownian motors. Europhys Lett 43:248-254.

\harvarditem{Parrondo et al.}{2002}{modmain} Parrondo JMR, de Cisneros BJ (2002) Energetics of Brownian motors: a review. Appl Phys A 75: 179-191.

\harvarditem{Paxton et al.}{2006}{pax} Paxton WF, Sundararajan S, Mallouk TE, Sen A (2006) Chemical locomotion. Angew Chem Int Ed 45: 5420-5429.

\harvarditem{Perthame and Souganidis}{2009a}{peth} Perthame B, Souganidis PE (2009a) Asymmetric potentials and motor effect: a homogenization approach. Ann Inst Henri Poincar\'{e} Anal Non Lin\'{e}aire 26:2055-2071.

\harvarditem{Perthame and Souganidis}{2009b}{pet} Perthame B, Souganidis PE (2009b) Asymmetric potentials and motor effect: a large deviation approach.
Arch Rat Mech Anal 193:153-169.

\harvarditem{Perthame and Souganidis}{2011}{petfg} Perthame B, Souganidis PE (2011) A homogenization approach to flashing ratchets. Nonlinear Differ Equ Appl 18:45-58.

\harvarditem{Peskin et al.}{1994}{mod3} Peskin CS, Ermentrout B, Oster G (1994) The correlation ratchet: a novel mechanism for generating directed motion by ATP hydrolysis. In: Mow VC et al. (eds) Cell Mechanics and Cellular Engineering. Springer, New York, pp 479-489.

\harvarditem{Prost et al.}{1994}{mod1} Prost J, Chauwin JF, Peliti L, Ajdari A (1994) Asymmetric pumping of
particles. Phys Rev Lett 72:2652-2655.




\harvarditem{Reimann and H\"{a}nggi}{2002}{bm} Reimann P, H\"{a}nggi P (2002) Introduction to the physics of Brownian motors. Appl Phys A 75: 169-178.

\harvarditem{Spirin}{2009}{sp} Spirin AS (2009) The ribosome as a conveying thermal ratchet machine. J Biol Chem 284(32):21103-21119.

\harvarditem{Vorotnikov}{2011}{vorr} Vorotnikov D (2011) The flashing ratchet and unidirectional transport of matter.
Discrete Contin Dyn Syst Ser B 16:963-971.


\end{thebibliography}
\end{document}